\documentclass[a4paper,11pt]{amsart}
\usepackage{amssymb}
\usepackage{mathrsfs}
\usepackage{amsmath,amssymb,amsthm,latexsym,amscd,mathrsfs}
\usepackage{indentfirst}
\usepackage{stmaryrd}
 \DeclareMathOperator{\Tr}{Tr}

 \newcommand{\ROM}[1]{\mathrm{\uppercase\expandafter{\romannumeral#1}}}

 \newtheorem{thm}{Theorem}[section]

\theoremstyle{definition}  
\newtheorem{ack}{Acknowledgements}   
 

\title[Willmore submanifolds in the unit sphere via isoparametric functions]{\textbf{Willmore submanifolds in the unit sphere via isoparametric functions}}

\author[Y.Q. Xie]{Yuquan Xie}
\address{Department of Mathematics, Hangzhou Normal University, Zhejiang 310036,
China} \email{yuqxie@pku.edu.cn}
\thanks {The project is supported by the NSFC ( No.11326071 ).}

 \subjclass[2000]{ 53A30, 53C42.}
\date{}
\keywords{Willmore submanifold, isoparametric hypersurface, focal submanifold.}
\begin{document}

\maketitle
\begin{abstract}
This paper is a continuation of \cite{TY12} and \cite{QTY13}. We show that both focal submanifolds of each isoparametric hypersurface in the sphere
with six distinct principal curvatures are Willmore.
\end{abstract}

\section{\textbf{Introduction}}
Let $x:N^n\rightarrow S^{n+p}$ be an immersion from an $n$-dimensional compact submanifold to an $(n+p)$-dimensional unit sphere.
Then $N$ is called a Willmore submanifold in $S^{n+p}$ if it is an extremal submanifold of the Willmore functional ($cf.$ \cite{Wan98})
\[   W(x)=\int_N (S-n|H|^2)^{\frac{n}{2}}dv.   \]
Here $S$ is the square norm of the second fundamental form of $x$, and $H$ is the mean curvature vector field.
In \cite{GLW01} and \cite{PW88}, the authors gave an equivalent condition for $N^n$ to be Willmore. In particular, if $N$ is minimal with constant $S$, the criterion for Willmore
reduces to a simple equation (see \eqref{ec} in section 2).

It follows immediately from the criterion \eqref{ec} that all the Einstein manifolds minimally
immersed in the unit sphere are Willmore submanifolds. However, there exist examples of minimal Willmore
submanifolds which are not Einstein, for example, Cartan's minimal isoparametric hypersurfaces. In addition,
\cite{Li01} characterized all the isoparametric Willmore hypersurfaces in the unite sphere.
In 2012, Tang and Yan \cite{TY12} proved that one of the focal submanifolds of OT-FKM
type is Willmore. 
Qian, Tang and Yan \cite{QTY13} extended this result, they showed that
\begin{thm}[\cite{QTY13}]
  Both the focal submanifolds of every isoparametric hypersurface in $S^{n+1}$ with four
  distinct principal curvatures are Willmore.
\end{thm}
Furthermore, they completely determine the focal submanifolds which are Einstein except for one
case (for details, see \cite{QTY13}). Recently, Tang and Yan \cite{TY13b} showed that the focal submanifolds with $g=4$
are all $\mathcal{A}-$manifolds but rarely Ricci parallel, except possibly for the only unclassified case.
For the case $g=6$, Li and Yan \cite{LY14} proved that none of them are Ricci parallel.

In this paper, we will give a new proof of theorem 1.1. Moreover,
we establish our result as follows
\begin{thm}
  Both the focal submanifolds of each isoparametric hypersurface in the unit sphere with six distinct principal curvatures
  are Willmore.
\end{thm}

Recall that an isoparametric hypersurface $M$ in the sphere is one whose principal curvatures
and their multiplicities are fixed constants.
In virtue of  M\"unzner's work \cite{Mun80} \cite{Mun81},
the number $g$ of distinct principal curvatures must be $1, 2, 3, 4$ or $6$, and there are
at most two multiplicities $\{m_1, m_2\}$ of principal curvatures ($m_1=m_2$, if $g$ is odd).
The classification problem has been completed except for one case (see \cite{Tho00} and \cite{Cec08} for excellent surveys and \cite{CCJ07}, \cite{Imm08}, \cite{Chi13},\cite{Miy13}, \cite{TY13a}, \cite{TXY14} for recent progresses).

The isoparametric hypersurfaces with $g=1,2,3$ were classified by Cartan to be homogeneous
\cite{Car39a} \cite{Car39b}. Clearly, in these cases,  the focal submanifolds are Willmore. For $g=6$, Abresch \cite{Abr83} proved that $m_1=m_2=1$ or $2$.
Dorfmeister and Neher \cite{DN85} showed that the isoparametric hypersurface is homogeneous in
the first case and Miyaoka \cite{Miy13} showed the same result in the second case.

Based on all the results mentioned above, we obtain the following
\begin{thm}
  All the focal submanifolds of the isoparametric hypersurfaces in the unit sphere are Willmore submanifolds.
\end{thm}

\section{\textbf{Notation and preliminary results}}

Let $N^n$ be a minimal submanifold in the unit sphere $S^{n+p}$ with constant square norm $S$ of the second fundamental form.
We choose a local field of orthonormal frames $e_1,\cdots, e_{n+p}$ in $S^{n+p}$ such that,
restricted to $N$, the vectors $e_1,\cdots,e_n$ are tangent to $N$ and, consequently,
the remaining vectors $e_{n+1},\cdots,e_{n+p}$ are normal to $N$.
Throughout this paper we will adopt the following ranges of indices:
\[   1\leq i,j,\cdots \leq n, \qquad n+1\leq \alpha, \beta, \gamma \cdots \leq n+p, \qquad 1\leq A, B, C, \cdots \leq n+p,  \]
and we shall agree that repeated indices are summed over the respected ranges.
With respect to the frame field in $S^{n+p}$ chosen above, let $\theta_1,\cdots,\theta_{n+p}$ be the  field of dual frames.
Then the structure equations of $S^{n+p}$ are given by
  \begin{align}
    &d \theta_A = \sum \omega_{AB}\wedge\theta_B,  \qquad\omega_{AB}= -\omega_{BA}, \nonumber\\
    &d \omega_{AB}= \sum \omega_{AC}\wedge\omega_{CB}-\theta_A\wedge\theta_B.\nonumber
  \end{align}
We restrict these form to $N$, then $$\theta_{\alpha}=0.$$
Since $0=d \theta_{\alpha}=-\sum \omega_{i\alpha}\wedge\theta_i$, by Cartan's lemma we may write
\begin{equation}
  \omega_{i\alpha}=\sum h_{ij}^{\alpha}\theta_j,\qquad h_{ij}^{\alpha}=h_{ji}^{\alpha}.\nonumber
\end{equation}
From these formulas, we obtain
  \begin{align}
    &d \theta_i = \sum \omega_{ij}\wedge\theta_j,  \qquad\omega_{ij}= -\omega_{ji},\\
    &d \omega_{ij}= \sum \omega_{ik}\wedge\omega_{kj}-\Omega_{ij}, \qquad \Omega_{ij}=\frac{1}{2}\sum R_{klij}\theta_k\wedge\theta_l,\\
    & R_{ijkl}=\delta_{ik}\delta_{jl}-\delta_{il}\delta_{jk}+\sum \left(h_{ik}^{\alpha}h_{jl}^{\alpha}-h_{il}^{\alpha}h_{jk}^{\alpha}\right). \label{Gauss-eq}
  \end{align}
Here $\omega_{ij}$, $\Omega_{ij}$  are the connection form and curvature form of $N$, respectively.
The Ricci curvature of $N$ is then given by
\begin{align}
  R_{ij}& = \sum_k R_{ikjk}=(n-1)\delta_{ij}+ \sum_{k,\alpha} \left(h_{ij}^{\alpha}h_{kk}^{\alpha}-h_{ik}^{\alpha}h_{kj}^{\alpha}\right) \nonumber\\
        & =(n-1)\delta_{ij}-\sum_{k,\alpha}h_{ik}^{\alpha}h_{kj}^{\alpha}. \label{Ricci}
\end{align}
The last equation successes, since $N$ is minimal.

Recall that if $N$ is minimal with constant $S$, the equivalent condition for $N$ to be Willmore
is given by ($cf.$ \cite{GLW01},\cite{TY12})
\begin{equation}
  \sum_{i,j}R_{ij}h_{ij}^{\alpha}=0, \qquad \text{for any $\alpha= n+1,\cdots, n+p$.}  \label{ec}
\end{equation}
Combining with \eqref{Ricci}, we conclude that $N$ is a Willmore submanifold, if it satisfies
\begin{equation}
  \sum_{i,j,k,\beta}h_{ik}^{\beta}h_{kj}^{\beta}h_{ij}^{\alpha}=0, \qquad \text{for any $\alpha= n+1,\cdots, n+p$.}
\end{equation}
Let $A_{\alpha}$ be the shape operator along the unit normal vector $e_{\alpha}$, that is,
$\left<A_{\alpha}(e_i), e_j\right>=h_{ij}^{\alpha}. $
We also denote  by $A_{\alpha}$ the corresponding matrix with respect to the orthonormal basis $e_1,\cdots,e_n$.
Since $N$ is minimal, the trace of the shape operator is
\begin{equation}
  \Tr A_{\alpha}=0\qquad \text{for any $\alpha= n+1,\cdots, n+p$.} \label{min-A}
\end{equation}
Using the fact
\begin{equation}
  \Tr\{\sum_{\beta} A_{\beta}^2A_{\alpha} \}
  =\sum_{i,j,k,\beta}h_{ik}^{\beta}h_{kj}^{\beta}h_{ij}^{\alpha}, \nonumber
\end{equation}
we see that $N$ is a Willmore submanifold in $S^{n+p}$, if it satisfies
\begin{equation}
  \Tr\{\sum_{\beta} A_{\beta}^2A_{\alpha} \}=0,
  \qquad \text{for any $\alpha= n+1,\cdots, n+p$.} \label{Wec-A}
\end{equation}

\section{\textbf{Proof of the Theorem}}
Let $N^n$ be a focal submanifold of the isoparametric hypersurface in the sphere $S^{n+p}$.
It is well known that $N$ is a minimal submanifold in $S^{n+p}$ with constant square norm $S$.
\subsection{Case $g=4$}
In this subsection, we shall give a new proof of Theorem 1.1.

In this case, as we know, the principal curvatures of $N$ with respect to any unit normal vector are
$ 1, 0, -1$. Therefore, 
if $\xi=\sum t_{\alpha}e_{\alpha}$ is any unit normal vector,
and its shape operator is denoted by $A$, then $A^3=A$.

By a further discussion (for details, one can find it in the previous version (arXiv0402272v1, Page 19) of \cite{CCJ07}),
we have
\begin{align}
  & A_{\alpha}=A_{\alpha}^3,\quad \text{for all $\alpha$},\\
  & A_{\alpha}=A_{\beta}^2A_{\alpha}+A_{\beta}A_{\alpha}A_{\beta}+A_{\alpha}A_{\beta}^2, \quad \text{for all $\alpha\neq\beta$.}
\end{align}
Hence for any $\alpha$, $\Tr A_{\alpha}=\Tr A_{\alpha}^3$, and for any $\alpha\neq\beta$
\begin{align}
  \Tr A_{\alpha}& =\Tr \{A_{\beta}^2A_{\alpha}+A_{\beta}A_{\alpha}A_{\beta}+A_{\alpha}A_{\beta}^2\}
                =3\Tr\{A_{\beta}^2A_{\alpha}\}.
\end{align}
Combining with \eqref{min-A}, we have for any $\alpha$,
\begin{align*}
  \Tr\{\sum_{\beta} A_{\beta}^2A_{\alpha} \} &=\Tr\{A_{\alpha}^3+\sum_{\beta\neq \alpha} A_{\beta}^2A_{\alpha} \}\\
        &=\Tr A_{\alpha}+\sum_{\beta\neq \alpha} \Tr\{A_{\beta}^2A_{\alpha}\}\\
        &=\Tr A_{\alpha}+\sum_{\beta\neq \alpha} \frac{1}{3}\Tr A_{\alpha}\\
        &=0,
\end{align*}
which yields that $N$ is Willmore. This completes the proof of theorem 1.1.
\subsection{Case $g=6$}
As mentioned before, $m_1=m_2=1$ or $2$ and the isoparametric hypersurfaces in these two cases are homogeneous.
Denote by $M_1, M_2$ the corresponding focal submanifolds.

For $m_1=m_2=1$, given $p\in M_1^5\subset S^7$, with respect to a suitable tangent orthonormal basis $e_1,\cdots,e_5$
of $T_pM_1$, Miyaoka \cite{Miy93} showed that the shape operators of $M_1$ are given by
\begin{equation*}
  \begin{matrix}
    A_6=\begin{pmatrix}
      \sqrt{3}  &   0   &   0   &   0   &   0\\
      0   &   \frac{1}{\sqrt{3}}   &   0   &   0   &   0\\
      0  &   0   &   0   &   0   &   0\\
      0   &   0   &   0   &   -\frac{1}{\sqrt{3}}   &   0\\
      0  &   0   &   0   &   0   &   -\sqrt{3}
    \end{pmatrix},\qquad    &
    A_7=\begin{pmatrix}
      0  &   0   &   0   &   0   &   \sqrt{3}\\
      0   &  0   &   0   &   \frac{1}{\sqrt{3}}   &   0\\
      0  &   0   &   0   &   0   &   0\\
      0   &   \frac{1}{\sqrt{3}}   &   0   &   0   &   0\\
      \sqrt{3}  &   0   &   0   &   0   &   0
    \end{pmatrix}.
  \end{matrix}
\end{equation*}
A direct computation leads to
\begin{equation*}
    A_6^2+A_7^2=\begin{pmatrix}
      6  &   0   &   0   &   0   &   0\\
      0   &   \frac{2}{3}   &   0   &   0   &   0\\
      0  &   0   &   0   &   0   &   0\\
      0   &   0   &   0   &   \frac{2}{3}   &   0\\
      0  &   0   &   0   &   0   &   6
    \end{pmatrix}.
\end{equation*}
It is not difficult to check that
\[  \Tr\{(A_6^2+A_7^2)A_{\alpha}\}=0,\qquad \alpha=6 \text{~or~} 7.    \]
Hence $M_1$ is a Willmore submanifold in $S^7$.

Similarly, for the focal submanifold $M_2^5$, the shape operators are given by ($c.f.$ \cite{Miy93} )
\begin{equation*}
  \begin{matrix}
    A_6=\begin{pmatrix}
      \sqrt{3}  &   0   &   0   &   0   &   0\\
      0   &   \frac{1}{\sqrt{3}}   &   0   &   0   &   0\\
      0  &   0   &   0   &   0   &   0\\
      0   &   0   &   0   &   -\frac{1}{\sqrt{3}}   &   0\\
      0  &   0   &   0   &   0   &   -\sqrt{3}
    \end{pmatrix},\qquad    &
    A_7=\begin{pmatrix}
      0  &   1   &   0   &   0   &   0\\
      1   &  0   &   0   &   -\frac{2}{\sqrt{3}}   &   0\\
      0  &   0   &   0   &   0   &   0\\
      0   &   -\frac{2}{\sqrt{3}}   &   0   &   0   &   1\\
      0  &   0   &   0   &   1   &   0
    \end{pmatrix}.
  \end{matrix}
\end{equation*}
Consequently,
\begin{equation*}
  \begin{matrix}
    A_6^2+A_7^2=\begin{pmatrix}
      4  &   0   &   0   &   -\frac{2}{\sqrt{3}}   &   0\\
      0   &  \frac{8}{3}   &   0   &   0   &   -\frac{2}{\sqrt{3}}\\
      0  &   0   &   0   &   0   &   0\\
      -\frac{2}{\sqrt{3}}   &   0   &   0   &   \frac{8}{3}   &   0\\
      0  &   -\frac{2}{\sqrt{3}}   &   0   &   0   &   4
    \end{pmatrix},\qquad    &
    \Tr\{(A_6^2+A_7^2)A_{\alpha}\}=0,\quad \alpha=6 \text{~or~} 7,
  \end{matrix}
\end{equation*}
which implies that $M_2$ is Willmore in $S^7$.

For the case $m_1=m_2=2$, the focal submanifolds $M_1^{10}, M_2^{10}$ are also homogeneous in $S^{13}$.
As asserted by Miyaoka \cite{Miy13}, the shape operators $A_{11}, A_{12}, A_{13}$ of $M_1$ are expressed
respectively by diagonal matrix
\begin{equation*}
    \begin{pmatrix}
      \sqrt{3}I  &   0   &   0   &   0   &   0\\
      0   &   \frac{1}{\sqrt{3}}I   &   0   &   0   &   0\\
      0  &   0   &   0   &   0   &   0\\
      0   &   0   &   0   &   -\frac{1}{\sqrt{3}}I   &   0\\
      0  &   0   &   0   &   0   &   -\sqrt{3}I
    \end{pmatrix},
\end{equation*}
and
\begin{equation*}
  \begin{matrix}
    \begin{pmatrix}
      0  &   0   &   0   &   0   &   \sqrt{3}J\\
      0   &  0   &   0   &   \frac{1}{\sqrt{3}}J   &   0\\
      0  &   0   &   0   &   0   &   0\\
      0   &   -\frac{1}{\sqrt{3}}J   &   0   &   0   &   0\\
      -\sqrt{3}J  &   0   &   0   &   0   &   0
    \end{pmatrix} ,   &
    \begin{pmatrix}
      0  &   0   &   0   &   0   &   \sqrt{3}I\\
      0   &  0   &   0   &   \frac{1}{\sqrt{3}}I   &   0\\
      0  &   0   &   0   &   0   &   0\\
      0   &   \frac{1}{\sqrt{3}}I   &   0   &   0   &   0\\
      \sqrt{3}I  &   0   &   0   &   0   &   0
    \end{pmatrix},
  \end{matrix}
\end{equation*}
where $I=\begin{pmatrix}
  1 & 0\\
  0 & 1
\end{pmatrix},
J=\begin{pmatrix}
  0 & -1\\
  1 & 0
\end{pmatrix}$.

Hence
\begin{equation*}
  \begin{matrix}
    A_{11}^2+A_{12}^2+A_{13}^2=\begin{pmatrix}
      9I  &   0   &   0   &   0   &   0\\
      0   &  I   &   0   &   0   &   0\\
      0  &   0   &   0   &   0   &   0\\
      0   &   0   &   0   &   I   &   0\\
      0  &   0   &   0   &   0   &   9I
    \end{pmatrix},\\
    \Tr\{(A_{11}^2+A_{12}^2+A_{13}^2)A_{\alpha}\}=0,\quad \alpha=11,12 \text{~or~} 13,
  \end{matrix}
\end{equation*}
which yields $M_1$ is Willmore in $S^{13}$.

For $M_2$, the shaper operator $A_{11}$ is the same as in $M_1$, and $A_{12}, A_{13}$ are given respectively by,
\begin{equation*}
  \begin{matrix}
    \begin{pmatrix}
      0  &   J   &   0   &   0   &   0\\
      -J   &  0   &   0   &   -\frac{2}{\sqrt{3}}J   &   0\\
      0  &   0   &   0   &   0   &   0\\
      0   &   \frac{2}{\sqrt{3}}J   &   0   &   0   &   J\\
      0  &   0   &   0   &   -J   &   0
    \end{pmatrix} ,   &
    \begin{pmatrix}
      0  &   -I   &   0   &   0   &   0\\
      -I   &  0   &   0   &   \frac{2}{\sqrt{3}}I   &   0\\
      0  &   0   &   0   &   0   &   0\\
      0   &   \frac{2}{\sqrt{3}}I   &   0   &   0   &  -I\\
      0  &   0   &   0   &   -I   &   0
    \end{pmatrix}.
  \end{matrix}
\end{equation*}
Thus,
\begin{equation*}
  \begin{matrix}
    A_{11}^2+A_{12}^2+A_{13}^2=\begin{pmatrix}
      5I  &   0   &   0   &   0   &   0\\
      0   &  5I   &   0   &   0   &   0\\
      0  &   0   &   0   &   0   &   0\\
      0   &   0   &   0   &   5I   &   0\\
      0  &   0   &   0   &   0   &   5I
    \end{pmatrix},\\
    \Tr\{(A_{11}^2+A_{12}^2+A_{13}^2)A_{\alpha}\}=0,\quad \alpha=11,12 \text{~or~} 13.
  \end{matrix}
\end{equation*}
This implies that $M_2$ is Willmore in $S^{13}$.

In summery, we conclude the theorem 1.2.


\begin{ack}
The author would like to express his hearty thanks to Professor Zizhou Tang for valuable discussions and constant encouragement,
and Dr. Wenjiao Yan and Chao Qian for their useful suggestions during the preparation of this paper.

\end{ack}

\providecommand{\bysame}{\leavevmode\hbox to3em{\hrulefill}\thinspace}
\providecommand{\MR}{\relax\ifhmode\unskip\space\fi MR }
\providecommand{\MRhref}[2]{%
  \href{http://www.ams.org/mathscinet-getitem?mr=#1}{#2}
}
\providecommand{\href}[2]{#2}

\end{document}